\documentclass[12pt]{article}
\usepackage{amssymb,amsmath,amsthm}
\usepackage{bm}
\usepackage{indentfirst}
\usepackage{amsfonts}
\usepackage{authblk}
\usepackage{times}
\usepackage{sectsty}
\usepackage{titlesec}
\providecommand{\keywords}[1]{\textbf{\textit{Keywords:}} #1}
\catcode240=13 \def ð{\u{g}} \catcode231=13 \def ç{\c{c}} \catcode246=13 \def ö{\"{o}} \catcode254=13 \def þ{\c{s}}
\catcode252=13 \def ü{\"{u}} \catcode253=13 \def ý{{\i}} \catcode221=13 \def Ý{\.{I}} \catcode199=13 \def Ç{\c{C}}
\catcode208=13 \def Ð{\u{G}} \catcode214=13 \def Ö{\"{O}} \catcode222=13 \def Þ{\c{S}} \catcode220=13 \def Ü{\"{U}}
\usepackage[a4paper,left=3cm,right=3cm,top=3cm,bottom=3cm]{geometry}
\usepackage{setspace}
\theoremstyle{definition}

\numberwithin{equation}{section}

\newtheorem{teo}{Theorem}

\newtheorem{cor}{Corollary}

\newtheorem*{pr}{Proof}
\newtheorem{ex}{Example}

\begin{document}

\title{A generalized integral transform and an alternative technique for solving linear ordinary differential equations}

\author[1 \footnote{Correspondence: ndernek@marmara.edu.tr}]{Nese Dernek}
\author[2]{Fatih Aylikci}

\affil[1]{Department of Mathematics, Marmara University, Istanbul, Turkey}
\affil[2]{Department of Mathematical Engineering, Yildiz Technical University, Istanbul, Turkey}

\maketitle

\begin{abstract}
In the present paper authors introduce the $\mathcal{L}_n$-integral transform and the $\mathcal{L}_n^{-1}$ inverse integral transform for ~$n=2^k$~,~$k \in \mathbb{N}$~, as a generalization of the classical Laplace transform and the $\mathcal{L}^{-1}$ inverse Laplace transform, respectively. Applicabi-\\lity of this transforms in solving linear ordinary differential equations is analyzed. Some illustrative examples are also given.
\end{abstract}


\section{Introduction, definitions and preliminaries}
The following Laplace-type the $\mathcal{L}_2$ transform
\begin{equation}
\mathcal{L}_2\{f(x);y\}=\int\limits_0^\infty x \exp(-x^2 y^2) f(x) dx
\end{equation}
was introduced by Yurekli and Sadek \cite{osman1}. After ~then Aghili, Ansari and Sedghi \cite{aghili} derived a complex inversion formula as follows
\begin{equation}
\mathcal{L}_2^{-1}\{\mathcal{L}_2\{f(x);y\}\}=\frac{1}{2\pi i}\int\limits_{c-i\infty}^{c+i\infty} 2 \mathcal{L}_2\{f(x);\sqrt{y}\}\exp(y x^2) dy
\end{equation}
where $\mathcal{L}_2\{f(x);\sqrt{y}\}$ has a finite number of singularities in the left half plane $Re ~y \leqslant c$. \\
\indent In this article, we introduce a new generalized integral transform as follows
\begin{equation}
\mathcal{L}_n\{f(x);y\}=\int\limits_0^\infty x^{n-1}\exp(-x^n y^n) f(x)dx \label{definition1}
\end{equation}
where $n=2^k, ~~ k \in \mathbb{N}$.\\
\indent The $\mathcal{L}_n$-transform is related to the Laplace transform by means
\begin{equation}
\mathcal{L}_n\{f(x);y\}=\frac{1}{n} \mathcal{L}\{f(x^{\frac{1}{n}});y^n\},
\end{equation}
where the Laplace transform is defined by
\begin{equation}
\mathcal{L}\{f(x);y\} \int\limits_0^\infty \exp(-xy)f(x)dx.
\end{equation}
\indent First we shall give several examples of the  $\mathcal{L}_n$-transforms of some elementary and special functions.
\begin{ex}
We show that
\begin{equation}
\mathcal{L}_n\{1;y\}=\int\limits_0^\infty x^{n-1}\exp(-y^n x^n)dx =\frac{1}{ny^n}.
\end{equation}
\end{ex}
\begin{ex}
We show that
\begin{equation}
\mathcal{L}_n\{x^k;y\}=\frac{1}{ny^{n+k}}\Gamma(\frac{k}{n}+1)  \label{example2}
\end{equation}
where $k,n \in \mathbb{N}$ and $k > -n$.
\indent Applying the definition of the $\mathcal{L}_n$-transform, we have
\begin{equation}
\mathcal{L}_n\{x^k;y\}=\int\limits_0^\infty x^{k+n-1} \exp(-y^n x^n)dx ,\label{example2.}
\end{equation}
where $k \in \mathbb{N}$. \\
\indent The integral on the right hand side may be evaluated by changing the variable of the integration from $x$ to $u$ where $x^n y^n = u$, and using Gamma function's relation in (\ref{example2.}), we obtain
\begin{equation}
\mathcal{L}_n\{x^k;y\} = \frac{1}{ny^{n+k}} \int\limits_0^\infty u^{\frac{k}{n}+1-1} \exp(-u)du = \frac{1}{ny^{n+k}} \Gamma(\frac{k}{n}+1).
\end{equation}
\end{ex}
\begin{ex}
We show that
\begin{equation}
\mathcal{L}_n\{\cos(ax^n);y\}=\frac{y^n}{n(y^{2n}+a^{2})}. \label{example3}
\end{equation}
\indent Using the definition of the $\mathcal{L}_n$-transform and calculating the Taylor expansion of the $\cos$ function in (\ref{example3}) we get
\begin{equation}
\mathcal{L}_n\{\cos(ax^n);y\}=\sum\limits_{m=0}^\infty (-1)^m \frac{a^{2m}}{(2m)!} \mathcal{L}_n\{x^{2mn};y\},
\end{equation}
where from (\ref{example2}) we have the following relation
\begin{equation}
\mathcal{L}_n\{x^{2mn};y\}= \frac{2m + 1}{ny^{n+2mn}}
\end{equation}
and then we obtain (\ref{example3})
\begin{equation}
\mathcal{L}_n\{\cos(ax^n);y\} = \frac{1}{ny^n} \sum\limits_{m=0}^\infty (-1)^m \frac{a^{2m}}{y^{2mn}} = \frac{y^n}{n(y^{2n}+a^2)}.
\end{equation}
\end{ex}
\begin{ex}
We show that
\begin{equation}
\mathcal{L}_n\{\sin(ax^n);y\}=\frac{a}{n(y^{2n}+a^2)}. \label{example4}
\end{equation}
\indent Using the linearity of the $\mathcal{L}_n$-transform and calculating the Taylor expansion of the $\sin$ function in (\ref{example4}) we get
\begin{equation}
\mathcal{L}_n\{\sin(ax^n);y\}=\sum\limits_{m=0}^\infty (-1)^m \frac{a^{2m+1}}{(2m+1)!} \mathcal{L}_n\{x^{(2m+1)n};y\}.
\end{equation}
Using the following relation in (\ref{example2}) of Example 1.2,
\begin{equation}
\mathcal{L}_n\{x^{(2m+1)n};y\}=\frac{\Gamma(2m+2)}{ny^{2m+2n}},
\end{equation}
we have
\begin{equation}
\mathcal{L}_n\{\sin(ax^n);y\}=\frac{a}{ny^{2n}}\sum\limits_{m=0}^\infty \frac{(-a^2)^m}{y^{2mn}}=\frac{a}{n(y^{2n}+a^2)}.
\end{equation}
\end{ex}
\begin{ex}
We show that
\begin{equation}
\mathcal{L}_n\{\exp(-a^nx^n);y\}=\frac{1}{n(y^n + a^n)}
\end{equation}
where ~$0 < Re (a) < Re (y) $. \\
\indent Using the definition of the $\mathcal{L}_n$-transform and calculating the Taylor expansion of the exponential function we have
\begin{equation}
\mathcal{L}_n\{\exp(-a^nx^n);y\}=\sum\limits_{m=0}^\infty (-1)^m \frac{a^{mn}}{m!} \mathcal{L}_n\{x^{mn};y\} \label{example5.}.
\end{equation}
Using the value
\begin{equation}
\mathcal{L}_n\{x^{mn};y\}=\frac{\Gamma(m+1)}{ny^{n+mn}}
\end{equation}
on the right hand side of (\ref{example5.}) we get
\begin{equation}
\mathcal{L}_n\{\exp(-a^nx^n);y\}=\frac{1}{ny^n}\sum\limits_{m=0}^\infty (-1)^m \frac{a^{mn}}{y^{mn}} = \frac{1}{n(y^n + a^n)}.
\end{equation}
\end{ex}
\begin{ex}
We show that
\begin{equation}
\mathcal{L}_n\{J_0(2a^{\frac{n}{2}}x^{\frac{n}{2}});y\}=\frac{1}{ny^n} \exp(-\frac{a^n}{y^n})
\end{equation}
where the function $J_0$ is the Bessel function of the first kind of order zero. \\
\indent Using the following Taylor expansion of the function $J_0(x)$,
\begin{equation}
J_0(x)=\sum\limits_{m=0}^\infty \frac{(-1)^m }{(m!)^2}\Big(\frac{x}{2}\Big)^{2m},
\end{equation}
we obtain
\begin{equation}
\mathcal{L}_n\{J_0(2a^{\frac{n}{2}}x^{\frac{n}{2}});y\}=\sum\limits_{m=0}^\infty \frac{(-1)^m a^{mn}}{m! \Gamma(m+1)} \mathcal{L}_n\{x^{mn};y\}. \label{ex3}
\end{equation}
We know the $\mathcal{L}_n$ transform of $f(x)=x^{mn}$ as
\begin{equation}
\mathcal{L}_n\{x^{mn};y\} = \frac{\Gamma(m+1)}{ny^{n+mn}}. \label{example1.6.1}
\end{equation}
Substituting the relation \ref{example1.6.1} into (\ref{ex3}) we obtain
\begin{equation}
\mathcal{L}_n\{J_0(2a^{\frac{n}{2}} x^{\frac{n}{2}});y\}=\frac{1}{ny^n}\sum\limits_{m=0}^\infty \frac{(-\frac{a^n}{y^n})^m}{m!} = \frac{1}{ny^n} \exp(-\frac{a^n}{y^n}).
\end{equation}
\end{ex}
\begin{ex}
We show that
\begin{equation}
\mathcal{L}_n\{x^{\frac{nv}{2}} J_v(2a^{\frac{n}{2}} x^{\frac{n}{2}});y\}=\frac{1}{n} a^{\frac{nv}{2}} y^{-n(v+1)}\exp(-\frac{a^n}{y^n}) \label{example7}
\end{equation}
where $Re (a) > 0 , Re v > -1$. \\
\indent Using the following Taylor expansion of $J_v(x)$, which is the Bessel function of the first kind of order $v$,
\begin{equation}
J_v(x)=\sum\limits_{m=0}^\infty \frac{(-1)^m }{m! \Gamma(m+v+1)}\Big(\frac{x}{2}\Big)^{2m + v},
\end{equation}
we obtain
\begin{equation}
\mathcal{L}_n\{x^{\frac{nv}{2}} J_v(2a^{\frac{n}{2}} x^{\frac{n}{2}});y\}=\sum\limits_{m=0}^\infty \frac{(-1)^m a^{nm + \frac{nv}{2}}}{m! \Gamma(m+v+1)} \mathcal{L}_n\{x^{mn + nv};y\}. \label{ex4}
\end{equation}
We can calculate the $\mathcal{L}_n$ transform of $f(x)=x^{mn + nv}$ function as follows
\begin{equation}
\mathcal{L}_n\{x^{mn + nv};s\} = \frac{\Gamma(m+v+1)}{ny^{n+mn+nv}}. \label{example1.7.1}
\end{equation}
Substituting the relation \ref{example1.7.1} into equation (\ref{ex4}) we obtain the assertion (\ref{example7}) of Example 1.7,
\begin{equation}
\mathcal{L}_n\{x^{\frac{nv}{2}} J_v(2a^{\frac{n}{2}} x^{\frac{n}{2}});y\}= \frac{a^{\frac{nv}{2}}}{ny^{nv + n}}\sum\limits_{m=0}^\infty \frac{(-\frac{a^n}{y^n})^m}{m!}=\frac{1}{n} a^{\frac{nv}{2}} y^{-n(v+1)}\exp(-\frac{a^n}{y^n}).
\end{equation}
\end{ex}
\begin{ex}
We show that
\begin{equation}
\mathcal{L}_n\{erfc(\frac{1}{2} ~a^{\frac{n}{2}} x^{-\frac{n}{2}});y\}=\frac{1}{n} y^{-n} \exp(-a^{\frac{n}{2}} y^{\frac{n}{2}}) \label{example8}
\end{equation}
where $Re (a) > 0$. \\
\indent Using the definition of the complementary error function $erfc(x)$,
\begin{equation}
erfc(x)=\frac{2}{\sqrt{\pi}}\int\limits_x^\infty \exp(-u^2)du,
\end{equation}
we get
\begin{equation}
\mathcal{L}_n\{erfc(\frac{1}{2} ~a^{\frac{n}{2}} x^{-\frac{n}{2}});y\}=\frac{2}{\sqrt{\pi}} \int\limits_0^\infty x^{n-1} \exp(-y^n x^n) \int\limits_{\frac{a^{n/2}}{2 x^{n/2}}}^\infty \exp(-u^2) du dx.
\end{equation}
Changing the order of integration, we obtain
\begin{equation}
\mathcal{L}_n\{erfc(\frac{1}{2} ~a^{\frac{n}{2}} x^{-\frac{n}{2}});y\}=\frac{2}{\sqrt{\pi}} \int\limits_0^\infty \exp(-u^2) \int\limits_{\frac{a}{2^{2/n} u^{2/n}}}^\infty x^{n-1} \exp(-y^n x^n) dx du
\end{equation}
and using the relation $\frac{d}{dx}(\exp(-y^nx^n))=-ny^nx^{n-1}\exp(-y^nx^n)$ we have
\begin{equation}
\mathcal{L}_n\{erfc(\frac{1}{2} ~a^{\frac{n}{2}} x^{-\frac{n}{2}});y\}=\frac{2}{\sqrt{\pi} n y^n} \int\limits_0^\infty \exp(-u^2) \exp(-\frac{y^n a^n}{4u^2}) du.
\end{equation}
Changing the variable from $u$ to $x$ according to the transformation $u= \frac{a^{n/2}}{2x^{n/2}}$ we find that
\begin{equation}
\mathcal{L}_n\{erfc(\frac{1}{2} ~a^{\frac{n}{2}} x^{-\frac{n}{2}});y\}=\frac{a^{\frac{n}{2}}}{2\sqrt{\pi} y^n} \mathcal{L}_n\{x^{-\frac{3n}{2}}\exp(-\frac{a^n}{4x^n});y\}.
\end{equation}
Using the Taylor expansion of exponential function and the $\mathcal{L}_n$-transform of \\$f(x)=x^{-mn - \frac{3n}{2}}$ we obtain
\begin{equation}
\frac{a^{\frac{n}{2}}}{2\sqrt{\pi} y^n} \mathcal{L}_n\{x^{-\frac{3n}{2}}\exp(-\frac{a^n}{4x^n});y\}=\frac{a^{\frac{n}{2}}}{2\sqrt{\pi} y^n} \sum\limits_{m=0}^\infty \frac{(-1)^m a^{mn}}{m! 4^m} \mathcal{L}_n\{x^{-mn - \frac{3n}{2}};y\} \nonumber
\end{equation}
\begin{equation}
=\frac{a^{\frac{n}{2}}}{2n\sqrt{\pi} y^n} \sum\limits_{m=0}^\infty \frac{(-1)^m a^{mn}}{m! 4^m} \frac{\Gamma(-m-\frac{3}{2}+1)}{y^{-mn - \frac{n}{2}}}.
\end{equation}
From the following Euler's reflection formula,
\begin{equation}
\Gamma(z)\Gamma(1-z)=\frac{\pi}{\sin(\pi z)},
\end{equation}
we get
\begin{equation}
\frac{a^{\frac{n}{2}}}{2n\sqrt{\pi} y^n} \sum\limits_{m=0}^\infty \frac{(-1)^m a^{mn}}{m! 4^m} \frac{\Gamma(-m-\frac{3}{2}+1)}{y^{-mn - \frac{n}{2}}}=\frac{a^{\frac{n}{2}} \pi}{2n\sqrt{\pi} y^n} \sum\limits_{m=0}^\infty \frac{(-1)^{2m+1} a^{mn} y^{mn + \frac{n}{2}}}{\Gamma(m+1) \Gamma(m+1+\frac{1}{2}) 4^m}
\end{equation}
and using the following duplication formula for Gamma function
\begin{equation}
\Gamma(z)\Gamma(z+\frac{1}{2})=2^{1-2z} \sqrt{\pi} ~\Gamma(2z)
\end{equation}
we obtain
\begin{equation}
\frac{a^{\frac{n}{2}} \pi}{2n\sqrt{\pi} y^n} \sum\limits_{m=0}^\infty \frac{(-1)^{2m+1} a^{mn} y^{mn + \frac{n}{2}}}{\Gamma(m+1) \Gamma(m+1+\frac{1}{2}) 4^m}=\frac{1}{n} y^{-n} \exp(-a^{\frac{n}{2}} y^{\frac{n}{2}}). \label{example8.1}
\end{equation}
Thus the assertion (\ref{example8}) follows from (\ref{example8.1}).
\end{ex}
\begin{ex}
We show that
\begin{equation}
\mathcal{L}_n\{erf(a^{\frac{n}{2}} x^{\frac{n}{2}});y\}=\frac{a^{\frac{n}{2}}}{n} y^{-n} (y^n + a^n)^{-\frac{1}{2}} \label{example9}
\end{equation}
where $-Re (a) < y , Re (y)>0$. \\
\indent Using the definition of $\mathcal{L}_n$-transform and the error function we have
\begin{equation}
\mathcal{L}_n\{erf(a^{\frac{n}{2}} x^{\frac{n}{2}});y\}=\frac{2}{\sqrt{\pi}}\int\limits_0^\infty x^{n-1} \exp(-y^n x^n) \int\limits_0^{a^{n/2} x^{n/2}} \exp(-u^2) du dx.
\end{equation}
Changing the order of integration and evaluating the inner integral we get
\begin{equation}
\mathcal{L}_n\{erf(a^{\frac{n}{2}} x^{\frac{n}{2}});y\}=\frac{2}{\sqrt{\pi}} \int\limits_0^\infty \exp(-u^2) \int\limits_{\frac{u^{2/n}}{a}}^\infty x^{n-1} \exp(-y^n x^n) dx du \nonumber
\end{equation}
\begin{equation}
=\frac{2}{\sqrt{\pi} n y^n} \int\limits_0^\infty \exp(-u^2(1+\frac{y^n}{a^n})) du.
\end{equation}
Changing the variable from $u$ to $x$ according to transformation $u\sqrt{1+\frac{y^n}{a^n}}=x$, ~we obtain the assertion (\ref{example9}),
\begin{equation}
\frac{2}{\sqrt{\pi} n y^n} \int\limits_0^\infty \exp(-u^2(1+\frac{y^n}{a^n})) du=\frac{a^{\frac{n}{2}}}{n} y^{-n} (y^n + a^n)^{-\frac{1}{2}}.
\end{equation}
\end{ex}
\begin{ex}
We show that
\begin{equation}
\mathcal{L}_n\{\exp(-ax^{2n});y\}=\frac{\sqrt{\pi}}{2n\sqrt{a}}~ \exp(\frac{y^{2n}}{4a}) erfc(\frac{y^{2n}}{2\sqrt{a}})
\end{equation}
provided that $Re ~a > 0$. \\
\indent Using the definition of the $\mathcal{L}_n$-transform we get
\begin{equation}
\mathcal{L}_n\{\exp(-ax^{2n});y\} = \int\limits_0^\infty x^{n-1} \exp(-y^nx^n - ax^{2n})dx. \label{example6.}
\end{equation}
Writing on the right hand side of (\ref{example6.})
\begin{equation}
-y^nx^n-ax^{2n}=-a(x^n + \frac{y^n}{2a})^2 + \frac{y^{2n}}{4a}
\end{equation}
and changing the variable
\begin{equation}
a^{1/2}(x^n + \frac{y^n}{2a})=u,
\end{equation}
using the definition of the complementary error function as follows, we deduce the assertion (\ref{example6.}),
\begin{equation}
\mathcal{L}_n\{\exp(-ax^{2n});y\}=\frac{\sqrt{\pi}}{2n\sqrt{a}}~\exp(\frac{y^{2n}}{4a}) erfc(\frac{y^n}{2\sqrt{a}}).
\end{equation}
\end{ex}
\begin{cor}
From the definition of the $\mathcal{L}_n$-transform the following identity hold true:
\begin{equation}
\mathcal{L}_n\{\exp(-ax^n)f(x);y\}=\mathcal{L}_n\{f(x);(y^n+a)^\frac{1}{n}\}
\end{equation}
where $Re ~a > 0$.
\end{cor}
\indent We now introduce a new derivative operator for the $\mathcal{L}_n$-transform and apply the operator to solve following ordinary differential equations:
\begin{equation}
xz'' - (2v+n-3)z' + x^{n-1}z=0~~ , ~~ n=2^k,k \in \mathbb{N},v>n,v=2^m +1, m \in \mathbb{N}\label{equation1}
\end{equation}
\begin{equation}
xz'' - (n^2 - 1)z' + x^{n-1}z=0~~,~~ n=2^k,~k=0,1,2,...  \label{equation2}
\end{equation}
\section{Some properties of the $\mathcal{L}_n$-transform}
In this section we will give some properties of the $\mathcal{L}_n$-transform that will be used to solve the ordinary differential equations (\ref{equation1})-(\ref{equation2}) given above. \\
\indent Firstly, we introduce a differential operator $\overline{\delta}$ (see \cite{osman4,osman5}) that we call the $\overline{\delta}$- derivative and define as
\begin{equation}
\overline{\delta}_x=\frac{1}{x^{n-1}}\frac{d}{dx} \label{oper1} , ~~ n=2^k,~k \in \mathbb{N}
\end{equation}
we note that
\begin{equation}
\overline{\delta}_x^2=\overline{\delta}_x \overline{\delta}_x=\frac{1}{x^{n-1}}\frac{d}{dx}\Big(\frac{1}{x^{n-1}}\frac{d}{dx}\Big)
=\frac{1}{x^{2n-2}}\frac{d^2}{dx^2} - \frac{(n-1)}{x^{2n-1}}\frac{d}{dx}.  \label{oper2}
\end{equation}
\indent The $\overline{\delta}$ derivative operator can be successively applied in a similar fashion for any positive integer power. \\
\indent Here we will derive a relation between the $\mathcal{L}_n$-transform of the $\overline{\delta}$-derivative of a function and the $\mathcal{L}_n$-transform of the function itself. \\
\indent Suppose that $f(x)$ is a continuous function with a piecewise continuous derivative $f'(x)$ on the interval $ [0,\infty) $. Also suppose that $f$ and $f'$ are of exponential order $\exp(\alpha^nx^n)$ as
$ x \rightarrow \infty $ where $\alpha$ is a constant. By using the definitions of $\mathcal{L}_n$-transform and the $\overline{\delta}$ derivative and integration by parts, we obtain
\begin{equation}
\mathcal{L}_n\{\overline{\delta}_x f(x);y\} = \int\limits_0^\infty \exp(-y^nx^n) f'(x)dx,
\end{equation}
\begin{equation}
\int\limits_0^\infty \exp(-y^nx^n) f'(x)dx=\lim\limits_{b \rightarrow \infty} f(x) \exp(-y^n x^n)|_0^b + ny^n \int\limits_0^\infty x^{n-1} \exp(-y^nx^n) f(x)dx.
\end{equation}
Since $f$ is of exponential order $\exp(\alpha^n x^n)$ as $x \rightarrow \infty$, it follows that
\begin{equation}
\lim\limits_{x \rightarrow \infty} \exp(-y^nx^n)f(x) = 0
\end{equation}
and consequently,
\begin{equation}
\mathcal{L}_n\{\overline{\delta}_x f(x);y\}=ny^n \mathcal{L}_n\{f(x);y\} - f(0^+). \label{operator1}
\end{equation}
Similarly, if $f$ and $f'$ are continuous functions with a piecewise continuous derivative $f''$ on the interval $ [0,\infty) $, and if all three functions are of exponential order $\exp(\alpha^n x^n)$ as $x \rightarrow \infty$ we can use (\ref{operator1}) to obtain
\begin{equation}
\mathcal{L}_n\{\overline{\delta}_x^2 f(x);y\}=n^2 y^{2n} \mathcal{L}_n\{f(x);y\} - n y^n f(0^+) - \overline{\delta}_x f(0^+). \label{operator2}
\end{equation}
Using (\ref{operator1}) and (\ref{operator2}) we get
\begin{equation}
\mathcal{L}_n\{\overline{\delta}_x^3 f(x);y\}=n^3 y^{3n} \mathcal{L}_n\{f(x);y\} - n^2 y^{2n} f(0^+) - ny^n \overline{\delta}_x f(0^+) - \overline{\delta}_x^2 f(0^+). \label{operator3}
\end{equation}
By repeated application of (\ref{operator1}) and (\ref{operator3}) we obtain the following theorem.
\begin{teo}
If $f,f',...,f^{(k-1)}$ are all continuous functions with a piecewise continuous derivative $f^{(k)}$ on the interval $[0,\infty)$, and if all functions are of exponential order $\exp(\alpha^n x^n)$ as $x \rightarrow \infty$ for some constant $\alpha$ then
\begin{equation}
\mathcal{L}_n\{\overline{\delta}_x^k f(x);y\} = (ny^n)^k \mathcal{L}_n\{f(x);y\} - (ny^n)^{k-1} f(0^+) - (ny^n)^{k-2} \overline{\delta}_x f(0^+) -
\end{equation}
\begin{equation}
... - ny^n \overline{\delta}_x^{k-2} f(0^+) - \overline{\delta}_x^{k-1} f(0^+)
\end{equation}
for $k \geq 1$, ~$k$ is a positive integer.
\end{teo}
The $\mathcal{L}_n$-transform defined in (\ref{definition1}) is an analytic function in the half plane $Re ~y> \alpha$. Therefore, $\mathcal{L}_n\{f(x);y\}$ has derivatives of all orders and the derivatives can be formally obtained by differentiating (\ref{definition1}). Applying the $\overline{\delta}$ with respect to the variable $y$ we obtain
\begin{equation}
\overline{\delta}_y \mathcal{L}_n\{f(x);y\}=\frac{1}{y^{n-1}} \frac{d}{dy} \int\limits_0^\infty x^{n-1} \exp(-y^nx^n) f(x) dx \nonumber
\end{equation}
\begin{equation}
=\frac{1}{y^{n-1}}\int\limits_0^\infty x^{n-1} (-x^n n y^{n-1} \exp(-y^n x^n)) f(x)dx = -n \mathcal{L}_n\{x^n f(x);y\}.
\end{equation}
If we keep taking the $\overline{\delta}$-derivative of (\ref{definition1}) with respect to the variable $y$, then we deduce
\begin{equation}
\overline{\delta}_y^k \mathcal{L}_n\{f(x);y\}=\int\limits_0^\infty x^{n-1} \overline{\delta}_y^k \exp(-y^n x^n) f(x) dx
\end{equation}
for $k \in \mathbb{N}$.
\begin{equation}
\int\limits_0^\infty x^{n-1} \overline{\delta}_y^k \exp(-y^n x^n) f(x) dx = \int\limits_0^\infty x^{n-1} \overline{\delta}_y^{k-1} [(-n) x^n \exp(-y^n x^n)]f(x)dx \nonumber
\end{equation}
\begin{equation}
=\int\limits_0^\infty x^{n-1} \overline{\delta}_y^{k-2}[(-n)^2 x^{2n} \exp(-y^n x^n)] f(x)dx \nonumber
\end{equation}
\begin{equation}
\ldots \nonumber
\end{equation}
\begin{equation}
=\int\limits_0^\infty x^{n-1} [(-n)^k x^{kn} \exp(-y^n x^n)] f(x)dx = (-n)^k \mathcal{L}_n\{x^{kn} f(x);y\}.
\end{equation}
As a result we obtain the following theorem.
\begin{teo}
If $f$ is piecewise continuous on the interval $[0,\infty)$ and is of exponential order $\exp(\alpha^n x^n)$ as $x \rightarrow \infty$, then
\begin{equation}
\mathcal{L}_n\{x^{kn} f(x);y\}= \frac{(-1)^k}{(n)^k} ~\overline{\delta}_y^k \mathcal{L}_n\{f(x);y\}
\end{equation}
for $k \geq 1$, $k$ is a positive integer.
\end{teo}
\begin{teo}
Let $\mathcal{L}_n\{f(x);y^{1/n}\}$, $n=2^k$, $k=0,1,2,...$ be an analytic function of $y$ except at singular points each of which lies to the left of the vertical line $Re~y=a$ and they are finite numbers. Suppose that  ~$y=0$ is not a branch point and $\lim\limits_{y \rightarrow \infty} \mathcal{L}_n\{f(x);y^{1/n}\}=0$ in the left plane $Re ~y \leq a$ then, the following identity
\begin{equation}
\mathcal{L}_n^{-1}\{\mathcal{L}_n\{f(x);y\}\}=\frac{1}{2\pi i} \int\limits_{a-i\infty}^{a+i\infty} n \mathcal{L}_n\{f(x);y^{1/n}\} \exp(yx^n) dy \nonumber
\end{equation}
\begin{equation}
=\sum\limits_{k=1}^m [Res\{n \mathcal{L}_n\{f(x);y^{1/n}\} \exp(yx^n);y=y_k\}] \label{theorem3}
\end{equation}
hold true for $m$ singular points.  \\
\pr
We take a vertical closed semi-circle as contour of integration. Using residues theorem and boundedness of $\mathcal{L}_n\{f(x);y^{1/n}\}$, we show that the identity (\ref{theorem3}) of Theorem 2.3 is valid. When $y=0$ is a branch point we take key-hole contour instead of simple vertical semi-circle. \\
\indent We assume that $\mathcal{L}_n\{f(x),y^{1/n}\}$ has a finite number of singularities in the left half plane $Re y \leqslant a$. Let $\gamma=\gamma_1 + \gamma_2$ be the closed contour consisting of the vertical line segment $\gamma_1$, which is defined from $a-iR$ to $a+iR$ and vertical semi-circle $\gamma_2$, that is defined as $|y-a|=R$. Let $\gamma_2$ lie to the left of vertical line $\gamma_1$. The radius $R$ can be taken large enough so that $\gamma$ encloses all the singularities of the $\mathcal{L}_n\{f(x);y^{1/n}\}$. Hence, by the residues theorem we have
\begin{equation}
\frac{1}{2\pi i} \int\limits_{a-i\infty}^{a+i\infty} n \mathcal{L}_n\{f(x);y^{1/n}\} \exp(yx^n) dy \nonumber
\end{equation}
\begin{equation}
=\frac{1}{2\pi i} \int\limits_{\gamma_1} n \mathcal{L}_n\{f(x);y^{1/n}\} \exp(yx^n) dy - \frac{1}{2\pi i} \int\limits_{\gamma_2} n \mathcal{L}_n\{f(x);y^{1/n}\} \exp(yx^n) dy \nonumber
\end{equation}
\begin{equation}
=\sum\limits_{k=1}^m [Res\{n \mathcal{L}_n\{f(x);y^{1/n}\} \exp(yx^n);y=y_k\}] - \frac{1}{2\pi i} \int\limits_{\gamma_2} n \mathcal{L}_n\{f(x);y^{1/n}\} \exp(yx^n) dy \label{proof3}
\end{equation}
where $y_1,y_2,\ldots,y_m$ are all the singularities of $\mathcal{L}_n\{f(x);y^{1/n}\}$. Taking the limit from both sides of the relation (\ref{proof3}) as $R$ tends to $+\infty$, because of the Jordan's Lemma, the second integral in the right tends to zero. \\
\indent Even $\mathcal{L}_n\{f(x);y^{1/n}\}$ has one branch point at $y=0$, we can use the identity (\ref{theorem3}). The proof of the proposition is similar to the proof of the Main Theorem in the paper \cite{aghili}, where we take $n=2^k,k \in \mathbb{N}$ instead of $n=2$. \\
\indent If the number of singularities is infinite, we take the semi-circles $\gamma_m$ which is centered at point $a$, with radius $R_m=\pi^2 m^2,m \in \mathbb{N}$.
\end{teo}
We illustrate the above Theorem by showing that the following examples.
\begin{ex}
We show that
\begin{equation}
\mathcal{L}_n^{-1}\{\frac{1}{y^{2n} + a^{2n}};x\}=\frac{n}{a^n} \sin(a^nx^n) \label{example11}
\end{equation}
where $Re ~a>0$. \\
\indent Using the assertion (\ref{theorem3}) of Theorem 2.3 we obtain
\begin{equation}
\mathcal{L}_n^{-1}\{\frac{1}{y^{2n}+a^{2n}};x\}=\sum\limits_{k=1}^2 Res[n \frac{1}{y^2 + a^{2n}} \exp(yx^n);y=y_k]  \label{example11.1}
\end{equation}
where the singular points are $y_k=\mp ia^n~,~k=1,2$ and then we have
\begin{equation}
Res[n \frac{1}{y^2 + a^{2n}} \exp(yx^n);ia^n] = n \lim\limits_{y \rightarrow ia^n} (y - ia^n) \frac{\exp(yx^n)}{y^2 + a^{2n}} = n \frac{\exp(ia^nx^n)}{2ia^n} \label{rezidu1}
\end{equation}
and similarly we have
\begin{equation}
Res[n \frac{1}{y^2 + a^{2n}} \exp(yx^n); - ia^n]=-n \frac{\exp(-ia^nx^n)}{2ia^n}. \label{rezidu2}
\end{equation}
Using the relations (\ref{rezidu1}) and (\ref{rezidu2}) we find the formula (\ref{example11}) from (\ref{example11.1}) as follows
\begin{equation}
\mathcal{L}_n^{-1}\{\frac{1}{y^{2n}+a^{2n}};x\} = \frac{n}{a^n} \frac{\exp(ia^nx^n)-\exp(-ia^nx^n)}{2i} \nonumber
\end{equation}
\begin{equation}
=\frac{n}{a^n} \sin(a^nx^n).
\end{equation}
\end{ex}
\begin{ex}
We show that
\begin{equation}
\mathcal{L}_n^{-1}\{\frac{1}{y^n} \exp(-\frac{a^n}{y^n});x\}=nJ_0(2a^{n/2} x^{n/2}) \label{example12}
\end{equation}
where $J_0$ is the Bessel function of order zero. \\
\indent Using the assertion (\ref{theorem3}) of Theorem 2.3 we have
\begin{equation}
\mathcal{L}_n^{-1} \{\frac{1}{y^n} \exp(-\frac{a^n}{y^n});x\}=Res[n \frac{1}{y} \exp(-\frac{a^n}{y}) \exp(yx^n),y=y_k]. \label{example12.1}
\end{equation}
From the following Taylor expansions of the exponential functions in (\ref{example12.1}),
\begin{equation}
n \frac{1}{y} \exp(-\frac{a^n}{y}) \exp(yx^n) = \frac{n}{y} \sum\limits_{m=0}^\infty (-1)^m \frac{a^{mn}}{m! y^m} \sum\limits_{k=0}^\infty \frac{y^k x^{nk}}{k!} \nonumber
\end{equation}
\begin{equation}
=\frac{n}{y}[1-\frac{a^n}{1! y} + \frac{a^{2n}}{2! y^2} -\frac{a^{3n}}{3! y^3}+...][1+ \frac{x^n y}{1!} + \frac{x^{2n} y^2}{2!} + \frac{x^{3n}}{3!}+...]
\end{equation}
we find that $Res[n \frac{1}{y} \exp(\frac{a^n}{y}) \exp(yx^n)]$ as the coefficient of the term $\frac{1}{y}$ as follows
\begin{equation}
Res[n \frac{1}{y} \exp(\frac{a^n}{y}) \exp(yx^n)=n [1-\frac{a^nx^n}{(1!)^2} + \frac{a^{2n} x^{2n}}{(2!)^2} - \frac{a^{3n} x^{3n}}{(3!)^2}+...] \nonumber
\end{equation}
\begin{equation}
=n\sum\limits_{m=0}^\infty (-1)^m \frac{(ax)^{mn}}{(n!)^2}=n J_0(2a^{n/2}x^{n/2}). \label{example12.2}
\end{equation}

Thus, we obtain from (\ref{example12.2}) and the formula (\ref{example12.1}), the assertion (\ref{example12}) of Example 2.12.
\end{ex}
\section{Application of the $\mathcal{L}_n$-transform to ordinary differential equations}
First we consider the ordinary differential equation (\ref{equation1})
for $v>n $~~and~~ $ \\ v=2^m + 1 ~~ , m=0,1,2,... $ \\
\indent Dividing ~(\ref{equation1}) by $x^{n-1}$, adding and subtracting the term $\frac{n-1}{x^{n-1}}z'$ we obtain
\begin{equation}
x^n\Big(\frac{1}{x^{2n-2}} z'' - \frac{n-1}{x^{2n-1}}z'\Big) + \frac{n-1}{x^{n-1}}z' - \frac{2v+n-3}{x^{n-1}}z' + z =0. \label{equation3}
\end{equation}
Using the definition of the $\overline{\delta}$-derivative given in (\ref{oper1}) and (\ref{oper2}), we can express (\ref{equation3}) as
\begin{equation}
x^n \overline{\delta}_x^2 z(x) - 2(v-1)\overline{\delta}_x z(x) + z(x)=0. \label{equation4}
\end{equation}
Applying the $\mathcal{L}_n$-transform to (\ref{equation4}) we find
\begin{equation}
\mathcal{L}_n\{x^n \overline{\delta}_x^2 z;y\} - 2(v-1)\mathcal{L}_n\{\overline{\delta}_x z ;y\} + \mathcal{L}_n\{z(x);y\}=0. \label{equation5}
\end{equation}
Using Theorem 2.1 for $k=1$ and $k=2$ in (\ref{equation5}) and performing necessary calculations we obtain
\begin{equation}
-\frac{1}{n} \overline{\delta}_y \mathcal{L}_n\{\overline{\delta}_x^2 z;y\} - 2(v-1)\mathcal{L}_n\{\overline{\delta}_x z;y\} + \mathcal{L}_n\{z;y\}=0,
\end{equation}
\begin{equation}
-\frac{1}{n}\frac{1}{y^{n-1}}\frac{d}{dy}(n^2 y^{2n} \overline{z}(y) - n y^n z(0^+) - \overline{\delta}_x z(0^+)) - 2(v-1)(n y^n \overline{z}(y) - z(0^+)) + \overline{z}(y)=0
\end{equation}
where $\overline{z}(y)=\mathcal{L}_n\{z(x);y\}$. We assume that $z(0^+)=0$. Thus, we obtain the following first order differential equation:
\begin{equation}
\overline{z}'(y) + (2(n+v-1) \frac{1}{y} - \frac{1}{ny^{n+1}})\overline{z}(y)=0. \label{equation6}
\end{equation}
Solving the first order differential equation (\ref{equation6}) we have
\begin{equation}
\overline{z}(y)=C \sum\limits_{m=0}^\infty (-1)^m \frac{1}{m! n^{2m} y^{mn + 2n + 2v - 2}}.
\end{equation}
Applying the $\mathcal{L}_n^{-1}$ transform we obtain
\begin{equation}
z(x)=C \sum\limits_{m=0}^\infty (-1)^m \frac{x^{mn + n + 2v -2}}{m! \Gamma(m + \frac{n+2v-2}{n} + 1) n^{2m-1}}
\end{equation}
where we use the following relations
\begin{equation}
\mathcal{L}_n\{x^k;y\}=\frac{\Gamma(\frac{k}{n}+1)}{ny^{n+k}} ~,~k = mn + n + 2v - 2
\end{equation}
and
\begin{equation}
\mathcal{L}_n^{-1}\{\frac{1}{y^{mn + n + 2v - 2 + n}}\}=\frac{n x^{mn + n + 2v -2}}{\Gamma(m+ 1 + \frac{2v-2}{n} + 1)}.
\end{equation}
Setting ~$\alpha=\frac{2v+n-2}{n}$~,~$ C=n^{-\frac{2v-2}{n}-2}$ ~ we obtain the solution of the ordinary differential equation (\ref{equation1})
\begin{equation}
z(x)=x^{\frac{n\alpha}{2}} J_{\alpha}(\frac{2}{n} x^{\frac{n}{2}})
\end{equation}
where $\alpha \in \mathbb{Z}$ because of the inequality $v > n ~(v,n \in \mathbb{N})$ and $J_{\alpha}$ is the Bessel function of the first kind of order $\alpha$.

\indent In the second step we will use the $\mathcal{L}_n$-transform for solving (\ref{equation2}).\\
\indent Dividing (\ref{equation2}) by $x^{n-1}$, adding and subtracting the term $\frac{n-1}{x^{n-1}}z'$ we obtain
\begin{equation}
x^n\Big(\frac{1}{x^{2n-2}z''(x) - \frac{n-1}{x^{2n-1}}z'(x)}\Big) + \frac{n-1}{x^{n-1}}z'(x) - (n^2-1) \frac{1}{x^{n-1}} z'(x) + z(x) =0.  \label{equation7}
\end{equation}
Using the definition of the $\overline{\delta}_x$-derivative (\ref{oper1}) and (\ref{oper2}) we can express (\ref{equation7}) as
\begin{equation}
x^n \overline{\delta}_x^2 z(x) - n(n-1)\overline{\delta}_x z(x) + z(x)=0. \label{equation10}
\end{equation}
Considering the following relations
\begin{equation}
\mathcal{L}_n\{x^n \overline{\delta}_x^2 z(x);y\}=-\frac{1}{n} \overline{\delta}_y \mathcal{L}_n\{\overline{\delta}_x^2 z(x);y\}=-2n^2 y^n \overline{z}(y) - ny^{n+1} \overline{z}'(y) + nz(0^+), \label{equation8}
\end{equation}
\begin{equation}
n(n-1)\mathcal{L}_n\{\overline{\delta}_x z(x);y\}= n(n-1) (ny^n \overline{z}(y) - z(0^+)) = n^2 (n-1) y^n \overline{z}(y) - n(n-1) z(0^+) \label{equation9}
\end{equation}
and applying the $\mathcal{L}_n$-transform to (\ref{equation10}) we obtain
\begin{equation}
\mathcal{L}_n\{x^n \overline{\delta}_x^2 z(x);y\} - n(n-1)\mathcal{L}_n\{\overline{\delta}_x z(x);y\} + \mathcal{L}_n\{z(x);y\}=0
\end{equation}
\begin{equation}
ny^{n+1} \overline{z}'(y) + [n^2(n+1) y^n - 1]\overline{z}(y) - n^2 z(0^+)=0 \label{equation12}
\end{equation}
where $\overline{z}(y)=\mathcal{L}_n\{z(x);y\}$.\\
\indent We may assume
\begin{equation}
z(0^+)=0.  \label{equation11}
\end{equation}
Solving the first order differential equation after substituting (\ref{equation11}) into (\ref{equation12}) we get
\begin{equation}
\overline{z}(y)=C y^{-n^2 - n} \exp(-\frac{1}{n^2y^n}). \label{equation13}
\end{equation}
Calculating the Taylor expansion of the exponential function in (\ref{equation13}) we have
\begin{equation}
\overline{z}(y)=C \sum\limits_{m=0}^\infty \frac{(-1)^m}{m! n^{2m}} \frac{1}{y^{n + nm + n^2}}. \label{equation14}
\end{equation}
Using the following relation
\begin{equation}
\mathcal{L}_n^{-1}\{\frac{1}{y^{n+nm+n^2}}\} = \frac{n x^{nm+n^2}}{\Gamma(m+n+1)}
\end{equation}
and applying the $\mathcal{L}_n^{-1}$ transform to (\ref{equation14}) we find
\begin{equation}
z(x)=C n^{n+1} x^{\frac{n^2}{2}} \sum\limits_{m=0}^\infty (-1)^m \frac{1}{m! \Gamma(m+n+1)} \Big( \frac{2x^{n/2}}{2n}\Big)^{2m+n}. \label{equation15}
\end{equation}
Setting $C=n^{-n-1}$ in (\ref{equation15}) we obtain the solution of the equation (\ref{equation2})
\begin{equation}
z(x)=x^{\frac{n^2}{2}} J_n\{\frac{2}{n}x^{\frac{n}{2}}\}
\end{equation}
where $J_n$ is the Bessel function of the first kind of order $n$.

\newpage

\end{document}